\documentclass[10pt]{amsart}
\usepackage{amsmath}
\usepackage{amsxtra}
\usepackage{amscd}
\usepackage{amsthm}
\usepackage{amsfonts}
\usepackage{amssymb}
\usepackage{eucal}

\newtheorem{lem}[subsection]{Lemma}

\theoremstyle{definition}

\theoremstyle{remark}

\newcommand{\thmref}[1]{Theorem~\ref{#1}}

\newcommand{\lemref}[1]{Lemma~\ref{#1}}
\newcommand{\defref}[1]{Definition~\ref{#1}}
\newcommand{\propref}[1]{Proposition~\ref{#1}}
\newcommand{\corref}[1]{Corollary~\ref{#1}}

\newcommand{\exref}[1]{Example~\ref{#1}}

\newcommand{\nc}{\newcommand}
\nc{\renc}{\renewcommand}
\nc{\ssec}{\subsection}
\nc{\sssec}{\subsubsection}
\nc{\on}{\operatorname}
\nc{\remm}[1]{\<{remark} \ \lbl{#1} \>{remark}}

\renc\k{\mathbf{k}}
\nc\ol{\overline} 
\nc\wh{\widehat}
\nc\tboxtimes{\wt{\boxtimes}}

 \nc{\equ}[1]{\underset{ {#1}}{\sim}}
 \def\cG{\check{\fG}}
\nc{\Aa}{{\mathbb{A}}}
 \nc{\Gg}{{\mathbb{G}}}

\nc{\Hh}{{\mathbb{H}}}

 \nc{\Nn}{{\mathbb{N}}}
\nc{\Pp}{{\mathbb{P}}}
\nc{\Rr}{{\mathbb{R}}}
\nc{\BV}{{\mathbb{V}}}
\nc{\BW}{{\mathbb{W}}}
\nc{\Zz}{{\mathbb{Z}}}
\nc{\Qq}{{\mathbb{Q}}}
\nc{\Ss}{{\mathbb{S}}}
\nc{\Cc}{{\mathbb{C}}}
\nc{\Ff}{{\mathbb{F}}}
\nc{\Oo}{{\mathcal{O}}}
\nc{\Mm}{{\mathcal{M}}}
\nc{\Ll}{{\mathbb{L}}}
\nc{\Lll}{{\mathbb{L}^*}}
\nc{\eo}{{\mathbf{\tau}}}
\nc{\dU}{{\overset{\bullet}{\bigcup}}{}}
\nc{\du}{{\overset{.}{\cup}}{}}
\nc{\dual}[1]{{\overset{\vee}{#1}{}}}

\nc{\cM}{{\check{\mathcal M}}{}}
 \nc{\oM}{{\overset{\circ}{\mathcal M}}{}}

 \nc{\fB}{{\mathfrak{B}}}
\nc{\tT}{{\mathfrak{T}}} 
\nc{\tTt}{{\mathfrak{T}_{\triangle}}}

\nc{\fF}{{\mathcal{F}}}

\nc{\bb}{{\mathbf{b}}}
\nc{\bc}{{\mathbf{c}}}
\nc{\bd}{{\mathbf{d}}}
\nc{\be}{{\mathbf{e}}}
\nc{\bj}{{\mathbf{j}}}
\nc{\bn}{{\mathbf{n}}}

\nc{\bph}{{\mathbf{\phi}}}
\nc{\bp}{{\mathbf{p}}}
\nc{\bq}{{\mathbf{q}}}
\nc{\bF}{{\mathbf{F}}}
\nc{\bu}{{\mathbf{u}}}
\nc{\bv}{{\mathbf{v}}}
\nc{\bx}{{\mathbf{x}}}
\nc{\bh}{{\mathbf{h}}}
\nc{\bs}{{\mathbf{s}}}
\nc{\by}{{\mathbf{y}}}
\nc{\bw}{{\mathbf{w}}}
\nc{\bA}{{\mathbf{A}}}
\nc{\bK}{{\mathbf{K}}}
\nc{\bI}{{\mathbf{I}}}
\nc{\bB}{{\mathbf{B}}}
\nc{\bG}{{\mathbf{G}}}
\nc{\bC}{{\mathbf{C}}}
\nc{\bD}{{\mathbf{D}}}
\nc{\bP}{{\mathbf{P}}}
\nc{\bH}{{\mathbf{H}}}
\nc{\bM}{{\mathbf{M}}}
\nc{\bN}{{\mathbf{N}}}
\nc{\bV}{{\mathbf{V}}}
\nc{\bU}{{\mathbf{U}}}
\nc{\bL}{{\mathbf{L}}}
\nc{\bT}{{\mathbf{T}}}
\nc{\bW}{{\mathbf{W}}}
\nc{\bX}{{\mathbf{X}}}
\nc{\bY}{{\mathbf{Y}}}
\nc{\bZ}{{\mathbf{Z}}}
\nc{\bS}{{\mathbf{S}}}
\nc{\ba}{{\mathbf{a}}}

\nc{\sA}{{\mathsf{A}}}
\nc{\sB}{{\mathsf{B}}}
\nc{\sC}{{\mathsf{C}}} 
\nc{\sF}{{\mathsf{F}}}
\nc{\sG}{{\mathsf{G}}}
\nc{\sK}{{\mathsf{K}}}
\nc{\sM}{{\mathsf{M}}}
\nc{\sO}{{\mathsf{O}}}
\nc{\sQ}{{\mathsf{Q}}}
\nc{\sP}{{\mathsf{P}}}
\nc{\sZ}{{\mathsf{Z}}}
\nc{\sfp}{{\mathsf{p}}}
\nc{\sr}{{\mathsf{r}}}
\nc{\sg}{{\mathsf{g}}}
\nc{\sff}{{\mathsf{f}}}
\nc{\sfb}{{\mathsf{b}}}
\nc{\sfc}{{\mathsf{c}}} 

\nc{\tA}{{\widetilde{{A}}}}
\nc{\tD}{{\widetilde{{A}}}}
\nc{\tH}{{\widetilde{{A}}}}
\nc{\tB}{{\widetilde{{B}}}}
\nc{\tg}{{\widetilde{\mathfrak{g}}}}
\nc{\tG}{{\widetilde{G}}}
\nc{\TM}{{\widetilde{\mathbb{M}}}{}}
\nc{\tO}{{\widetilde{\mathsf{O}}}{}} 
\nc{\TZ}{{\tilde{Z}}}
\nc{\tx}{{\tilde{x}}}
\nc{\tf}{{\tilde{f}}}
\nc{\tz}{{\tilde{\zeta}}}
\nc{\tmu}{{\tilde{\mu}}}
\nc{\td}{{\tilde{d}}}
\nc{\tX}{{\widetilde{X}}}
  
   \nc{\E}{{\mathop{\operatorname{\rm E }}}}
 \nc{\Mor}{{\mathop{\operatorname{\rm Mor \,}}}}
\nc{\Ob}{{\mathop{\operatorname{\rm Ob \,}}}}
  \nc{\Sym}{{\mathop{\operatorname{\rm Sym}}}}
   \nc{\Aut}{{\mathop{\operatorname{\rm Aut}}}}
 \nc{\Spec}{{\mathop{\operatorname{\rm Spec}}}}
\nc{\Ker}{{\mathop{\operatorname{\rm Ker}}}}
 \nc{\dom}{{\mathop{\operatorname{\rm dom}}}}
\nc{\End}{{\mathop{\operatorname{\rm End}}}}
 \nc{\Hom}{\on{\Hom}} 
 \nc{\GL}{{\mathop{\operatorname{\rm GL}}}}
 \nc{\Id}{{\mathop{\operatorname{\rm Id}}}}
 \nc{\rk}{{\mathop{\operatorname{\rm rk}}}}
\nc{\irk}{{\mathop{\operatorname{\rm i-rk}}}}
 \nc{\length}{{\mathop{\operatorname{\rm length}}}}
\nc{\supp}{{\mathop{\operatorname{\rm supp}}}}
\nc{\val}{{\rm val}}

\nc{\valr}{{\rm val_{rv}}}
\nc{\valrv}{\valr}

\nc{\res}{{\mathop{\operatorname{\rm res}}}}

\nc{\ac}{{\mathop{\operatorname{\rm ac}}}}

\def\meet{\cap}
\def\union{\cup}

\def\si{\sigma}

\def\G{\Gamma}
\def\<{\begin}
 \def\>{\end}

\nc{\tV}{{\widetilde{{V}}}}
\nc{\hb}[1]{\hbox{#1}}
\def\rv{{\rm rv}}
 
\nc{\seq}[1]{\stackrel{#1}{\sim}}

\nc{\oeq}[1]{\underset{#1}{=}}

\def\inv {{^{-1}}}
 \def\beq#1{\begin{equation} \label{#1}}   
\def\eeq{\end{equation}}

  \def\Vv{\mathbb V}

\def\iso{\simeq}

\def\I-{{I^{-}}}  \def\Lm{{L^{-}}}   \def\Tm{{T^{-}}}   \def\M-{{M^{-}}}  \def\A-{{A^{-}}}
\def\prf{\begin{proof}}
\def\eprf{\end{proof} }

\def\acl{\mathop{\rm acl}\nolimits}
 \def\dcl{\mathop{\rm dcl}\nolimits}

\author{Ehud Hrushovski}

 \address{\newline Institute of Mathematics, the Hebrew
 University of Jerusalem, Givat Ram, Jerusalem, 91904, Israel.} 
 
 \email{ehud@math.huji.ac.il}

\usepackage{color}

\def\lbl#1{  \label{#1}  }

\nc{\Claim}[1]{{\noindent \bf Claim{ #1 }}}
  \nc{\pr}{{\mathop{\operatorname{\rm pr}}}}

\nc{\Mmm}{{(1+\Mm)}}

     \nc{\wt}{{\mathop{\operatorname{\rm wt}}}}

  \title{On finite imaginaries}

\begin{document}

 \maketitle

 \def\Vvs{ { \Vv ^s}}
 \def\tri{{\mathsf B}}
 \def\tres{{\mathbf{t}}}
 \def\tTs{\widetilde{T_\si}}
 \def\fG{{\mathcal G}}
 \def\fE{{\mathcal E}}
 \def\fC{{\mathcal C}}
\def\fD{{\mathcal D}}
\def\fGV{{\fG}_{\Vv}}
\def\ObG{{\rm Ob}_{\fG}}  
\def\ObGV{{\rm Ob}_{\fG_{\Vv}}}
\def\MorG{{\rm Mor}_{\fG}}
 %---------------------------------------------------------------------------------------
 \def\cB{{\mathfrak B}_{\rm cl}}
 \nc{\av}{{\mathop{\operatorname{\rm av}}}}
 \def\RV{   {\rm RV}   }
\nc{\K}{{ \rm K}}

\def\VF{{ \mathbb{\rm VF}}}
 \def\VFr{{  \VF \rv }}
\def\VFn{{ \mathbb{\rm \VF [n]}}}
\def\RVn{{ \mathbb{\rm \RV_{\T}[n]}}}
 \def\VFm{  \mu {\VF}}
 \def\RVm{  \mu {\RV}} 
 \def\VFnm{\mu \VFn}
\def\RVnm{{ \mu \VFm}}
 
     \def\RVmn{{\rm RV_{\mu}[n]}}
 
   \def\VFmn{\mu \VF [n]}
   \def\VFgn{\VFm_\G[n]}

   \def\Gm{{ \Gamma \mu}}

    \def\VFmb{ \VFm _{\rm bdd}}
 \def\VFmbg{ \VFm _{\G;\rm bdd}}
 \def\mVFr{ {\VFm} {\rv}}
 
 \def\RVmb {\RVm _{\rm bdd}}
 \nc{\JG}{{\mathop{\operatorname{{ Jcb_{\G}  }}}}}
 \nc{\JRV}{{\mathop{\operatorname{{ Jcb_{\RV}  }}}}}
 \nc{\JVF}{{\mathop{\operatorname{{ Jcb  }}}}}

\def\rvi{ \rv^{-1}} 
\def\fCg{{\fC}^{RV}_\G}
 \def\fCgb{{\fC}^{RV}_{\G^{bdd}}}
 \def\fCga{ {\G_A}}
 \def\fCgab{{\fC}^{RV}_{\G_A^{bdd}}}
  \def\emp{\emptyset}
  \renc\b{\beta}
 \renc\a{\alpha}
 
  \nc{\tF}{{\widetilde{F}}}
      \nc{\tU}{{\widetilde{U}}}
 \def\m{\setminus}
\def\bt{\beta}
\newtheorem{claim}[subsection]{Claim}

\nc{\conjrv}{\underset{rv}{\sim}}
\def\rvinv{{$\conjrv$-invariant }}

\def\rvi{ \rv^{-1}}
\def\A{\alpha} 
 
%----------------t4
\def\tchi{\widetilde{\chi}}
\def\tC{\widetilde{C}}
\def\d{\delta}
\def\wX{\overline{X}}
\def\wY{\overline{Y}}
 
\def\Gpos{{\G^{\geq 0}}}
  %---------------------------------------------------------------------------------------
 \def\Om{{\Omega}}
 \def\om{{\omega}}

  \def\AA{{\mathcal A}}
  \def\AAF{F \AA}
\def\IA{{I \AA}}
\def\AAk{ \Omega \AA_\k}
\def\IAk{I Fn(\k,\AG(\St))} 
\def\RR{{\mathcal R}}
\def\CC{{\mathcal C}}
 \newcommand{\isomto}{\overset{\rightarrow}{\cong}}

\def\VFrf{{{\VF}_{full}^{\RV}}}
\def\VFs{{{\VF}_{sp}^{RV}}}
\def\fnr{{Fn^{\RV}}}

\def\cRV{{{\RV}}}

\def\rh{{\vartheta}}
 \nc{\SG}{{\mathop{\operatorname{{  K_+ }}}}}
  \nc{\SGe}{{\mathop{\operatorname{{  K_+^{eff} }}}}}
 \def\SK{\SG}
 \def\SKe{\SGe}
 \nc{\AG}{{\mathop{\operatorname{ { K} }}}}
 \nc{\AGe}{{\mathop{\operatorname{ { K^{eff}} }}}}
 \def\AGM{\mu \AG}
 
 \def\Isp{{\rm I_{sp}}}
 \def\Ivf {{\rm I_{vf}}}
 \def\Ivfm{{\rm I_{vf_\mu}}}
 
 \def\intn{\int_{[n]}}
 \def\intns{\int_{[n]}^{\SG}}
  \def\ints{\int^{\SG}}
 
 \def\pRV{{\rm pRV }}
 \def\uX{{\underline{X}}}     
\def\L{{\Lambda}}
   \def\uL{{\underline{\L}}}       
 
      \def\uY{{\underline{Y}}}

    \def\lemm#1{    \begin{lem} \lbl{#1} }
\def\T{{\mathfrak t}}
\def\V{{\rm V}}

\def\Tdcl{{\dcl^T}}
\def\Ttp{{{\rm tp}^T}}

\def\tP{{\widetilde{P}}}
\def\tQ{{\widetilde{Q}}}
\def\tY{{\widetilde{Y}}}

\def\fCgb{{\rm \G^{bdd}_A}}
\def\fCgf{{\rm \G^{fin}_A}}

\def\St{{\rm RES}}
\def\sd{\ltimes}
\def\vol{{\rm vol}}

\def\RCF{{\rm RCF}}
\def\RCVF{{\rm RCVF}} 
\def\DOAG{{\rm DOAG}}
\def\ACVF{{\rm ACVF}}
\def\RES{{\rm RES}}
\def\VAL{{\rm VAL}}
\def\pCF{{\rm pCF}}

%revised March 17, 2007 following referee report.
 
\def\fnF{ {^{\subseteq F}}}   \def\SF{{S \fnF}} 
\<{abstract}
We study finite imaginaries in certain valued fields, and prove a conjecture of Cluckers and Denef.
 \>{abstract}

\maketitle

\<{section}{Introduction}
In their beautiful paper \cite{CD}, Cluckers and Denef
study actions of linear algebraic groups on 
varieties over local fields of high residue characteristic.  The orbits of these actions,
known to be finite in number, can be viewed as 
imaginary elements of the   theory $HF_0$ of Henselian fields of   residue characteristic zero.
Cluckers and Denef relate them to finite imaginaries of   a certain extension,
 $T_\infty^{(d)}$ (cf. \ref{tinfd}).
 They formulate
 a ``tameness'' conjecture (2.16) on the nature of such imaginaries, and show that it 
 is true for imaginaries arising from algebraic  group actions over number fields
 (\cite{CD} Theorem 1.1).  From this they obtain  consequences for orbital integrals (Theorem 1.2).   

Sections   3 of this note contain  a proof of  Conjecture 2.16 in general (\thmref{2.16}).  This in turn is a special case
of a more general description of imaginaries in this theory, and indeed in a wider class
of theories of Henselian fields.  However as Denef suggested finite imaginaries
have a certain autonomy that permits their direct classification.  The proof of the finite case
follows  the lines of the general (unpublished) proof, and has the merit of  containing most the main ideas while avoiding technicalities.  \S 2 contains some general observations on finite imaginaries, while \S 3 describes them for the theory  $T_\infty^{(d)}$.

The rest of the paper contains two further comments on \cite{CD}, from different points of view.   Consider  theories $T$  
of Henselian valued fields  of residue characteristic $0$.
The data of \cite{CD}  consists of an algebraic group  $G$ and a homogeneous space $V$ for 
$G$, in the sense of algebraic geometry.  This means that an action of $G$ on 
$V$ is given, and over an algebraically closed field the action is transitive.  
But for a given field $L \models T$ the action need not be transitive, and the  question concerns the space of orbits  of $G(L)$ on $V(L)$. The goal is to show that $V(L)/G(L)$ reduce to imaginaries of the residue
field and the value group; and indeed to special, ``tame'' imaginaries (see below.)

 The data $G,V$ is geometric
(i.e. quantifier-free) and group-theoretic.  Both of these qualities are lost in the reductions
of $G(F)/V(F)$   to the residual sorts given by one of the above methods.  
In \S 4 we describe another method that does not lose track of the group theory, 
and is geometric in the sense of being independent of a particular completion of the theory of Henselian fields.  We use the theory ACVF; but we cannot simply interpret   the implied quantifiers of $G$-conjugacy on $V$ in terms of the quantifier elimination of ACVF, 
  since the resulting quotient would be trivial.  Reformulating the question in terms of groupoids
  does allow us to work in ACVF without trivializing the problem.   We illustrate this  in the special   case   $V=G/T$, where $T$ is a torus.
 
In \cite{CD} the residue field is pseudo-finite, the value group  essentially $\Zz$-like, and 
 tameness is defined in concrete terms;  an imaginary is tame if in the Denef-Pas language it is definable over the sorts $\k$ (the residue field)
 and $\G / n\G$ (where $\G$ is the value group); but not $\G$ itself. 
% In particular, the elements of $\G$ are not tame. 
  In the final section we attempt to understand the role of tameness from a more abstract viewpoint.   In \cite{CDM}, Cherlin, Van den Dries and Macintyre
consider   imaginaries coding the Galois groups of finite   extensions of a field $F$.  
\footnote{cf.    \S 1 of \cite{CH2} for an account and further references.  In these references, 
 the field $F$ is PAC, but this condition is irrelevant here.}    In  the context of $HF_0$,
 we note that Galois imaginaries satisfy a strong (and more symmetric) form of tameness. 
 We also show that  $V(L)/G(L)$ is analyzable over the Galois sorts, and hence over the tame imaginaries.
 \footnote{It follows
 that in the Denef-Pas language they are definable over the sorts $\G/n\G$ and $\k$.
 Asides from this remark, we use the intrinsic valued field structure in this paper, 
 and do not split $\RV$.} 
   We leave open the question of whether analyzability can be replaced by
 internality.  In an appendix, we discuss Galois sorts for theories more general than fields.
 
\S1-3, \S 4, and \S 5  can each be read independently of the others, and of \cite{CD};
however \S 5   (like \cite{CD}) relies  on \cite{serre}.  
{\em  In this paper, by  ``definable'' we mean   0-definable, i.e. invisible parameters are not allowed. }

\>{section}

\<{section}{Finite sorts} 

We will consider first order theories $T$ in many-sorted languages $L$.   $T$ need not be complete
but everything we do can easily be reduced to the complete case. Terms like ``surjective'' applied to definable relations and functions mean: provably in $T$. 
Thus a definable surjection $D \to D'$ means a definable $R \subset D \times D'$,
such that in any model $M \models T$, $R(M)$ is the graph of a surjection $D(M) \to D'(M)$.

A sort $S$ of $L$ is called {\em finite} if $S(M)$ is finite for all $M \models T$.

A family ${\mathcal F}$ of sorts is said to be {\em closed under products } if  the product of two sorts in ${\mathcal F}$ is definably 
  isomorphic to a third.

\<{defn}  
If $S$ is a sort of $L$ and $\{S_j\}_{j \in J}$
is a family sorts, we say $S$ is {\em dominated} by  $\{S_j\}_{j \in J}$ if there
exists a definable 
$D \subseteq S_{j_1} \times \ldots \times S_{j_k}$ and a definable   surjective 
map $f: D \to S$, with $j_1,\ldots,j_k \in J$.  

\>{defn}

Equivalently, in any   $M \models T$, $S(M) \subseteq \dcl(\union_{j \in J} S_j(M))$.

%A family of sorts  $\{S_j\}_{j \in J}$ is called {\em dominant} if every sort is dominated by
%it.  Equivalently, in any $M \models T$, $\dcl(\union_{j \in J} S_j(M)) = M$.

In this language, Conjecture 2.16  can be stated as follows.
\<{defn}[\cite{CD}] \lbl{tinfd}  Let   $T_\infty^{(d)}$ be the theory of Henselian fields with pseudo-finite residue field of characterstic $0$, and (dense) value group elementarily equivalent to $$(\{\frac{a}{b}: a \in \Zz, b \in \Nn, (b,d)=1 \} \leq \Qq,<,+) $$   \>{defn}

\<{thm}\lbl{2.16} Every finite imaginary sort of $T_\infty^{(d)}$ is
dominated by the set  of   sorts consisting of the residue field $\k$ and the finite value group quotient sorts $\G/ n \G$. \>{thm}

This will be proved at the end of \S 3.

\<{defn} \lbl{twist}  Let $S,F$ be sorts of $T$ with $F$ finite.  Define the imaginary sort $\SF$
to be the sort of partial functions $F \to S$. \>{defn}\

Explanation:  by definition, for some $N \in \Nn$,  $T \models |F| \leq N$.  
Then ${\SF} = \union_{n \leq N} {\SF}_n$ where ${\SF}_n$ is the sort of partial functions
$F \to S$ with $n$ element domain.       By  identifying a function with its graph, 
an $n$-element set of pairs, and then identifying a set of pairs with a tuple of pairs
up to $Sym(n)$, we see that ${\SF}_n$ embeds naturally into  the imaginary sort  $(F \times S)^n / Sym(n)$.  Observe in particular that ${\SF}$ is dominated by $S,F$. 

We view an element of $\SF$ as a tuple of elements of $S$, indexed by a finite subset of $F$
in place of a finite subset of $\Nn$.  We thus refer to the sorts $\SF$ as ($Aut(F)$-) twisted Cartesian powers of $S$.
 
 \lemm{aux2}  Let $T$ be a theory.  Assume given 
 a family ${\mathcal S}$  of sorts, 
  and a family ${\mathcal F}$ of  finite sorts,
  closed under products. 
    Let $T'$ be obtained from $T$ by naming the elements of each $S \in {\mathcal F}$.
I.e. $T'=Th((M,c_j)_{j \in J})$ for some $M \models T$ and enumeration $(c_j : j \in J)$ of
 ${\mathcal F} = \union_{S \in {\mathcal F}} S(M)$. 
 If   $T'$ eliminates imaginaries to the sorts  ${\mathcal S}$, then
  $T$ eliminates imaginaries to the sorts  $\{{\SF}: S \in {\mathcal S}, F \in  {\mathcal F} \}$    \>{lem}

\prf Observe that $(S_1 \times S_2)  \fnF$   embeds naturally into to $S_1 \fnF \times {S_2} \fnF$.
Thus  we may assume  ${\mathcal S}$ closed under Cartesian products.
 
Let $M \models T$ and let $e$ be an imaginary element of $M$.   We have to find
$h \in {\SF}$ for some $S,F$ with $\dcl(e)=\dcl(h)$.  
 
 Let $M'$ be an expansion of $M$ to a model of $T'$.  By assumption, in $M'$,
 $\dcl_{M'}(e)=\dcl_{M'}(g)$ for some tuple $g$ from the sorts ${\mathcal S}$; so we may
 assume $g \in S(M)$, $S \in {\mathcal S}$.  It follows that there exists $ F \in {\mathcal F}$
 and $d \in F(M)$ such that $\dcl(e,d)=\dcl(g,d)$ in the sense of $M$.  So $g = H(e,d)$
 and $e=G(g,d)$ for some definable functions $G,H$.  We can restrict the domains of
 $H,G$ to any given definable set containing $(e,d)$ (respectively $(g,d)$.)  
 So we may assume that $G(H(x,y),y) = x$.  Let $H_e(y) = H(e,y)$.  So $H_e$
 is a function with nonempty domain contained in $F$, and range in $S$.   So $H_e$ is coded by some $e$- definable element $h$ of ${\SF}$.  On the other 
 hand $e$ is determined by $H_e$, in fact $e=G(H_e(y),y)$ for any $y \in \dom(H_e)$.
 So $\dcl(e)=\dcl(h)$.  \eprf

\<{lem} \lbl{rem-p}  Suppose ${\mathcal F}$ is a collection of finite 
cyclic groups; and the relation:  ``there exists a definable surjective homomorphism $A \to B$''
 is a directed partial ordering on ${\mathcal F}$, i.e any two  groups 
 in ${\mathcal F}$ are definable homomorphic images of  a third.   Then the condition  of being closed
    under products in \lemref{aux2} can be dispensed with. \>{lem}
    
 \prf    Let $B_1, \ldots, B_n \in {\mathcal F}$.  There exists  $B \in  {\mathcal F}$
 and definable  surjective maps
 $h_i: B \to B_i$.  For any  $S \in {\mathcal S}$, this gives rise to 
 a surjective $h: B^n \to \Pi_{i=1}^n B_i$ and hence 
 a definable injection
 $S^{B_1 \times \ldots \times B_n} \to S^{B^n}$, 
$y \mapsto  y \circ h $.   Thus a   twisted powers of a set $S$ by a product of   sorts of ${\mathcal F}$ embeds into   a twisted power $X^{B^k}$ by a power of single such $B$.   
    
 Now we use the fact that $B$ is cyclic, say of order $d$.  Then the map $B \times [1,\ldots,d]^k \to B^k$, $(b,(a_1,\ldots,a_k)) \mapsto (b^{a_1},\ldots,b^{a_k})$, is a surjective map.
Let $N = d^k$ and fix a bijection $\{1,\ldots,N\} \to  [1,\ldots,d]^k $.  Then 
as above $S^{B^k}$ embeds definably into $S^{B \times [1,\ldots,d]^k} = (S^B)^{N}$,
and similarly for partial functions.  So the sorts $S^B$, $B \in {\mathcal F}$ 
dominate the sorts $S^{B_1 \times \ldots \times B_n}$, $B_1,\ldots, B_n \in {\mathcal F}$.
\eprf     

 \<{lem}  \lbl{FV} Let $T_1,T_2$ be two theories, and let $T=T_1 \times T_2$ be their 
 Feferman-Vaught product, so that a model of $T$ is a product of a model of $T_1$ with
 one of $T_2$, and a relation is a Boolean combination of products of relations.  Assume $T_1,T_2$ eliminate imaginaries and that every sort of $T_2$ is linearly ordered.
 Then $T$ eliminates imaginaries. \>{lem}
 
\prf  Let $M \models T$ and let $R$ be an $M$- definable relation on $M$. 
Say $e$ is a canonical code for $R$, possibly imaginary.  
 Then $M=M_1 \times M_2$, $M_i \models T_i$, and $R$ is a finite disjoint union of products
 of a definable set of $M_1$ with a definable set of $M_2$.  
 Let $B_2$ be the Boolean algebra of definable subsets of $M_2$ generated
 by the sections $R(a) = \{y: (a,y) \in R\}$.  This algebra is finite; let $\{R^i_2: i=1,\ldots,k \}$ be an enumeration
 of the atoms (without repetitions.)  Let $e^i_2$ be a  canonical parameter for $R^i_2$.
 Then the set  $\{e^i_2\}$ is
$e$-definable.   Since we assumed each sort of $T_2$ is linearly ordered, each $e^i_2$ is actually $e$-definable.  Let $R^i_1 = \{x: \{x \} \times R^i_2 \subseteq R\}$.  Then $R^i_1$
is clearly $e$-definable, and $R= \union_i    R^i_1 \times R^i_2$.   Let $e^i_1$
be a canonical code for $R^i_1$.  Then $\dcl(e) = \dcl((e^i_1,e^i_2)_i)$. So $e$ is coded
by a real element.  \eprf

\ssec{Elimination of finite imaginaries}
   
Let $T$ be a theory in a many-sorted language $L$, with  sorts $S_i \  (i \in I)$.   Let
$\I- \subset I$, and 
let $\Lm$  be the language consisting of the sorts $S_i \ (i \in \I-)$ and the relations among them;
let $\Tm = T | \Lm$.    If $M \models T$,
$\M-$ will denote the restriction of $M$ to the sorts $S_i (i \in \I-)$.

 If $C \subseteq C' \subseteq M \models T$,
we say $C'$ is   stationary over $C$ if 
$\dcl^{eq}(C') \meet \acl^{eq}(C) = \dcl^{eq}(C)$, i.e. every imaginary element
that is definable over $C'$ and algebraic over $C$ is definable over $C$.   It is clear that
if $C''$ is stationary over $C'$, and $C'$ over $C$, then so is $C''$ over $C$.  
A  type $p=tp(c/C)$ is called {\em stationary} if $C \union \{c\}$ is stationary over $C$.  

\<{lem} \lbl{stat1}  

(1) A type $p$ over $C$ is stationary if and only if for any $C$-definable equivalence relation $E$ on a $C$-definable set $D$ with finitely many
classes, if $p(x) \models D(x)$ then $p$ chooses one of the classes of $E$, i.e. $p(x) \union p(y) \models xEy$. 

(2) If there exists a $C$-definable type $\bar{p}$ over $M$ extending $p$, then $p$
is stationary.

%(3) \lemm{stat1}  Let $C$ be a substructure of $M \models T_\infty^{(d)}$.  Let $B$ be a $C$-definable ball, with no proper definable sub-balls.  Then $B$ generates a complete type
%$p$ over $C$, and $p$ is stationary.
\>{lem}
 
\prf  (1)  If $p=tp(c/C)$ is not stationary, then for some $C$-definable function $f$ (possibly with imaginary
values), $f(c)$ is algebraic over $C$ but not definable.  Define $x E y \iff f(x)=f(y)$.
Then $E$ divides the solutions of $p$ to finitely many classes, but more than one.
Conversely if   a $C$-definable equivalence relation $E$ divides $p$ into finitely many classes
(but more than one), then $c/E$ is an imaginary element in $\acl^{eq}(C) \m \dcl^{eq}(C)$.

(2) Let $E$ be a  $C$-definable equivalence relation   with finitely many
classes.     Since $\bar{p}$ is
a complete type over $M$ there is a unique 
class $D$ of $E$ such that $(x E a) \in \bar{p}$ iff $a \in D$.   By definability of the type, there exists a formula $\theta(y)$ over $C$
such that $(xEa) \in \bar{p}$ iff $\theta(a)$.  Thus $D$ is $C$-definable, by $\theta$;
and $p \models D(x)$.   \eprf

\<{lem} \lbl{fei}  Assume 

(1) $T$ admits elimination of $S_i$ -quantifiers for each $i \in I \m \I-$.  $\Tm$ is stably embedded in $T$.

(2)  Let  $M \models T$ be a countable model.  
Then there exists $C$ containing $\M-$ and stationary over $\M-$,
such that $\acl(C) \prec M$.
 
(3)  For $A \leq M \models T$,   let $T_A = Th(M,a)_{a \in A}$.  If 
     $F$ is a finite $T_A$-definable set   then there exists a finite $T_A$-definable set $F'$
  of $\M-$   and a definable bijection $g: F \to F'$. 
   
Then every finite imaginary sort of $T$ is definably isomorphic to one of $\Tm$.   \>{lem}

\prf  Quantifier-elimination  will be used in the background, to avoid disagreement about the notion of a definable
set between $T$ and $\Tm$.

Let $S = D  / E $ be a finite imaginary sort, where $D$ is a definable set in $T$; let $\pi: D \to S$ be the canonical map.  We may assume all elements of $S$ realize the same type.
Let $M \models T$.  View $\M-$ as a subset of $M$, and  let  $C$ be as in (2),
and $N = \acl(C)$.  
Then by (2), $N \prec M$.  Since $S$ is finite,  $S(N)=S(M)$, so each class of $E$ has a representative
in $N$.  Thus there exists $e \in C$ and a finite $e$-definable set $H_e \subseteq D$,
meeting every $E$-class.   Using (3), let $W$ be a finite $e$-definable subset of $\M-$
and $h_e: W \to H_e$ a definable bijection.   
 
By stable embeddedness (1), $W$  
is actually defined over some $e' \in \M-$.  Write $W=H'_{e'}$. 
We have an induced $e$-definable surjection $\psi_{e'}:H'_{e'}  \to S$.  
But there are only
finitely many maps $H'_{e'} \to S$, hence all are algebraic over $\M-$; by the stationarity assumption (2),
since $\psi_{e'}$ is $C$-definable it is also $\M-$-definable.  By enlarging $e'$ we may assume
it is $e'$-definable.

Let $H'$ be a definable set  of $\Tm$, and 
  $\psi: H' \to D$  a definable  map, such that % $e' \in Q$,
 $H'_{e'} =\{x: (x,e') \in H' $, $\psi_{e'} (x)  = \psi  (x,e')$.  Then the composition $\pi \circ \psi: H' \to  S$
is surjective.    Let $E'$ be the kernel, i.e. define
  $E'(x,y) \iff \pi \circ \psi (x) = \pi \circ \psi(y)$.  Then we have
a definable bijection $H'/E' \to S$.   By stable embeddedness again, $E'$ is $\Tm$-definable
and $H'/E'$ is an imaginary sort of $\Tm$.   \eprf

\>{section}

\<{section}{Finite sorts of $T_\infty^{(d)}$}

 Let $PF$ be the theory of pseudo-finite fields, $PF_0 = PF + \hbox{ ``characteristic 0 "}.$
   Let $F^{Gal}_n$ be a finite imaginary sort
 whose elements code the  elements of the Galois
group of the unique field extension of order $n$.  Let 
${\mathcal F^{Gal}} = \{F_2,F_3,F_4,\ldots \}$. 
Let $PF'$ be the theory obtained from $PF$ by naming the elements of ${\mathcal F^{Gal}_n}$
for each $n$.  
 It was shown in  \cite{CH} that $PF'$
eliminates imaginaries.  Since   $ F^{Gal}_n$ surjects canonically onto $F^{Gal}_m$ when $m$ divides $n$,
 \lemref{rem-p} applies.   Hence, by \lemref{aux2}   we have:  

\<{example} \lbl{ex-finite}    $PF$ eliminates imaginaries to the level of ${\mathcal F^{Gal}_n}$-twisted powers of the field sort. \>{example}

 The finite group $\Zz/n \Zz$ admits elimination of imaginaries to the level
of subsets of $\Zz / n \Zz$.    To see this, it suffices to note that any subgroup  $H$ of the automorphism group $G$
 of $\Zz/n\Zz$ has the form $H = \{g \in G: gY=Y \}$ for some
$Y \subseteq  \Zz/n\Zz$.  Indeed, we have $G =  (\Zz / n  \Zz)^* $ so that 
$H \subseteq G \subset \Zz / n  \Zz$, and we can let $Y = H$.

\<{example} \lbl{group}  Let $T(d)$ be the theory of ordered Abelian groups, divisible by all primes
$q$ with $(q,d)=1$, and such that for $p | d$ we have $\G/p^m \G \iso \Zz/p^m$,
and moreover an isomorphism $\G/p^m \G \iso \Zz/p^m$ is given as part of the language
(say by a predicate for the pullback to $\G$ of $1 \in \Zz/p^m$; then each element of
$\G/p^m \G$ becomes definable, as a multiple of the distinguished generator.)
Then $T(d)$ admits elimination of imaginaries. 

Let $ \Zz^{(d)} = \{a/b \in \Qq: a,b \in \Zz, b \neq 0, (b,d)=1 \}$.  
 Then $Th(\Zz^{(d)})$ in the ordered
group language admits EI to the sorts $\G$  together with the $\G/n$ and
the sort of subsets of $\G/n$ (where $n$ can be taken to be a power of $d$.)  \>{example}

\proof The theory $T(d)$ admits elimination of quantifiers, and it is easy to classify the
definable subsets of $\G$ and the definable functions $\G \to \G$ (in one variable)
and see explicitly that they are coded.  This suffices in general, cf. \cite{acvf1}, and shows
that $T(d)$ eliminates imaginaries.  

It follows from \lemref{aux2} that   $Th(\Zz^{(d)})$ eliminates imaginaries
to the level of twisted products $\G^{\G/n_1 \times \ldots \times \G/n_2}$. 
This reduces to $\G ^{(\G/n)^k}$ and again, by \lemref{rem-p}, to  the
sorts $\G^{\G/n}$.  Now a function $\G/n \to \G$ carries the same information as
a subset of $\G$ of size $\leq n$ (the image), together with a partial ordering of a subset
$\G/n$.   As remarked above this reduces to subsets of $\G/n$.

\<{cor} \lbl{tame} Let $T$ be the model-theoretic disjoint sum of $PF_0$ and $Th(\Zz^{(d)})$:
$T$ has two sorts $\k,\G$, with relations $+,\cdot,0,1$ on $\k$, $+,<,0$ on $\G$; such that $\k \models PF_0$
and $\G \models Th(\Zz^{(d)})$.   Then $T$ eliminates imaginaries to the sorts $\G$ together
with the sorts $S(n_1,n_2)$ for $n_1,n_2 \in \Nn$, where $S(n_1,n_2) =   \k^{F_1 \times F_2}$, with $F_1  = F^{Gal}_{n_1}, F_2 = \G/n_2$.
   \>{cor}
   
\prf  Let ${\mathcal S}$ be the set of Cartesian products of these sorts.  

\Claim{}  ${\mathcal S}$
 is closed under twisted powers
by $\G/n$. 

\prf  Using the linear ordering, a function $\G/n \to \G$
can be coded by an $n$-tuple of elements of $\G$ together with a partial ordering on
$\G/n$, namely the pullback of the linear ordering on $\G$.  
 In turn this partial ordering can be coded by a function from $\G/n$ to a subset
of the prime field of $\k$.   This shows that functions $\G/n \to \G$ are coded in 
${\mathcal S}$.  On the other hand a function
$\G/n \to Y ^{\G/m}$ (with $Y= \k^{F_1}$) can be viewed as a function $\G/n \times \G/m \to Y$,
and handled using \lemref{rem-p}.  \eprf

 Thus by \lemref{aux2} it suffices to prove that the theory
$T'$ obtained by naming the elements of each $\G/n$ 
eliminates imaginaries to the sorts   $\G$  together
with $\k^{F^{Gal}_{n}}$.  But this follows  from \lemref{FV} together with \exref{ex-finite}
and the first part of  \exref{group}.
\eprf

\ssec{Finite imaginaries in $T_\infty^{(d)}$}.  

We take the theory $T_\infty^{(d)}$  with the sorts
$\VF,\k,\G$, as well as the  sorts $\G/2,\G/3,\ldots$ and all twisted product sorts $\k^{F^{Gal}_{n_1} \times \G/n_2}$.  The  sort of these last two kinds will be called $I_t$.   The language includes the field structure on $\VF$, 
a valuation map  $\VF \m \{0 \} \to \G$,
predicates for the valuation ring $\Oo$ and its maximal ideal $\Mm$, 
 the residue homomorphism $\Oo/\Mm \to k$, and finally 
  a group isomorphism $\VF^* / (1+ \Mm) \to (\k^* \times \G)$
splitting the canonical maps $k^* \to \VF^* / (1+ \Mm) \to \G$.  The projection
to $\G$ is thus canonical, while the projection to $\k^*$ is the   Denef-Pas ``angular component'' map.

  Let
 $\I- =\{\G\} \union I_{t}$, $I = \I- \union \{\VF\}$ where $\VF$ is the valued field sort.
We will now show that the hypotheses (1-3) of \lemref{fei}   are valid for
$T_\infty^{(d)}$.  
 
\lemref{fei} (1)   follows from  Denef-Pas elimination of quantifiers,  \cite{pas}.  
 Stable embeddedness is clear from the form of the  quantifier-elimination; see \cite{acvf1}
 for an identical proof in the case of  ACVF.

\lemm{1c}  Let $C \subseteq M \models T_\infty^{(d)}$.  Assume $\M- \subseteq C$
and the maps $\val, \ac$ restricted to $C$ are  onto the value group and residue field of $M$.  
Then $\acl(C) \prec M$.  \>{lem}

\prf  Let $N = \acl(C) \meet \VF(M)$.   Then $N$ is Henselian, so by the Ax-Kochen, Ershov principle for the Denef-Pas language applied to  $T_\infty^{(d)}$, we see that $N,\VF(M)$ are elementarily equivalent,
and by model completeness, $N \prec \VF_M$.  It follows immediately
from the surjectiveness that $\acl(C) \prec M$.  \eprf

Now to prove \lemref{fei} (2),   let   $M \models T=T_\infty^{(d)}$.  Let $\RV = \VF / (1 + \Mm)$, $\rv: \VF \to \RV$ the quotient map.  Let $\A- =  \{\rv(m): m \in M\}$.  From the point of view of $M$, $\rv(m)$ is equi-definable with the pair $(\val(m),\ac(m)$, and so $\A-$ can be identified with
$\val(M) \times \k^*(M) \subseteq \dcl(\M-)$.  But $\A-$ can also be considered as a substructure of a model of $\ACVF$, in the language
considered in \cite{HK}.  

Let $C$ be a maximal subset of $\VF(M)$ such that 
$C \union \M-$ is stationary over $\M-$.   Let $A$ be the definable closure of $\A- \union C$
in the sense of $\ACVF$.  

We have to show that $\acl(C  \union \M-) \prec M$, and by 
\lemref{1c} it suffices to show that $(\val ,\ac) | C$ is surjective.   In other words
given $m \in M$, to find $c \in C$ with $\val(c) = \val(m), \ac(c) =\ac(m)$.  Let   $\b = \rv(m)$,
  $B_\b = \rv \inv (\b)$.  We have to find $c \in B_\b(C):= B_\b \meet C$.    

\Claim{1}   
$B_\b$ is not transitive in $\ACVF_A$.  In other words some proper sub-ball of $B_\b$ is
$\ACVF_A$-definable.
\prf
Suppose for contradiction that $B_\b$ is not transitive in $\ACVF_A$.
In this case, for any polynomial $F$, $\rv \circ F$ is constant on $B_\b$.  (\cite{HK} Lemma 3.47.)
%\lemref{rv-transitive} = 3.47
So $\ac (F) $ and $\val(F)$ are constant.  By 
  Denef-Pas quantifier elimination, we see that $B_\b$ is also
transitive in $(T_\infty^{(d)})_{\{\M- \union C\}}$.  (Compare \cite{HK} Lemma 12.1)
Let $p$ be the unique type over $\M- \union C$ of elements of $B_\b$.  

Moreover there exists a $\b$-definable type $p_\b(x)$ over $M$ concentrating on elements of $B_\b$.
Namely, first pick a $\val(m)$-definable type $r(t)$ of elements of $\G$, concentrating on intervals
$(\val(m),\val(m)+\epsilon) \subseteq \G$.  Over $M$, the type $r$  can be

$$r(t) = \{t \in \G, t > \val(m) \} \union \{t< s: s \in \G(M), s > \val(m)\} \union 
\{(\exists t')(nt'=t): n=1,2,3,\ldots \}$$
Then let 
$$p_\b(x) = \{x \in B_\b \} \union \{ \val(x-u) \models r: u \in B_\b(M) \} \union \{\ac(x-u)=1: 
u \in B_\b(M) \} $$
It is clear that this is a complete, consistent type over $M$ and is $(\ac(m),\val(m))$-definable.

By  \lemref{stat1}  (2), $p = p_\b | \M-\union \{C\}$ is stationary.
Choose any $c \in B_\b(M)$.  Then $tp(c/\M-\union \{C\})=p$.
So $C \union \{c\} \union \M-$ is stationary over $\M-$, contradicting the maximality of $C$.  
\eprf

Hence $B_\b$ 
contains a proper $\ACVF_A$-definable closed ball.    In this case by \cite{HK} Lemma 3.39 
%\lemref{red1c1}
$B_\b$ contains an $\ACVF_A$-definable point $d$.  So  $d \in \dcl(\M-\union \{C\})$ and
hence $d \in C$.

This finishes the proof of (2).
  
(3) We may take $F \subset \VF^n$.  $T_\infty^{(d)}$ is algebraically bounded
in the sense of \cite{vddries}, so $F$
   is contained in a finite $ACVF_A$-quantifier-free definable set $F'$.
 This  reduces us to the same lemma for ACVF; for a proof,
see for example \cite{HK} Lemma 3.9.  

 \<{cor}  Every finite definable imaginary set of  $T_\infty^{(d)}$ can be definably
embedded into  some power of the twisted product sorts
 $\k^{F^{Gal}_{n_1} \times \G/n_2}$  \>{cor}

\prf      \lemref{fei} reduces this to imaginary sorts of $(T_\infty^{(d)})^-$; so 
by \corref{tame} it suffices to show this for a finite definable   
$D \subset \G^m \times P$, where $P$ is a power of twisted product sorts.   
We can use induction on the cardinality, so we may assume $D$ is not 
 the union of two proper definable subsets.  Since $\G$ is linearly ordered, it follows
 that the projection $D \to \G^m$ has a one-point image.    Thus
$D$ projects injectively to a product of twisted product sorts.   
\eprf

By the remark below \defref{twist}, the  twisted product sorts
 $\k^{{F^{Gal}_{n_1} \times \G/n_2}}$  are dominated by the   sorts 
 $\k$, $\G/2,\G/3,\ldots$, and ${F^{Gal}_{n}}$.  These Galois sorts
are dominated by $\k$; hence every finite definable imaginary set of 
 $T_\infty^{(d)}$ is dominated by the tame sorts $\k$ and $\G/2,\G/3,\ldots$.
 
  This proves \thmref{2.16}.
 
 \>{section}

\<{section}{Groupoids} %and reductions, and Galois sorts}%Third proof}

A {\em   groupoid} is a category $\G=(Ob_\G, Mor_\G)$ in which every morphism is invertible.  We will consider definable   groupoids with a single isomorphism type.  See \cite{cigha},
though the use of definable groupoids there is different.  A {\em morphism} between groupoids
is a quantifier-free definable functor.  

Given a groupoid $\G$ defined without quantifiers in  a theory $T$,
one obtains  an equivalence relation $E_\G$, defined uniformly over   $T=Th(L)$
for any definably closed subset $L$ of a model of $T$.  Namely 
the equivalence relation of $\G$-isomorphism on  the objects of $\G$:  for $a, b \in Ob_\G(L)$,
$a E_\G ^L b $ iff $Mor(a,b)(L) \neq \emptyset$.   Let $Iso(\G;L)$ be the quotient
$Ob_\G(L) / E_\G^L$, i.e. the set of isomorphism classes of $\G(L)$.  
\footnote{ This 
connection between groupoids and imaginaries is different from the one considered in \cite{cigha}.  The approach in this section is apparently in the spirit of  stacks.  }

Let $HF_0$ be the theory of   Henselian fields with   residue
field of characteristic $0$.
  This can be viewed as the theory of definably closed substructures of models of $ACVF_\Qq$ (with trivially valued $\Qq$.)   Fix $L \models HF_0$.  
We will use only quantifier-free formulas, and notions such as $\dcl$ will refer to $ACVF_L$.

We wish to reduce a given quantifier-free definable groupoid $\G$ over ACVF to a groupoid $\G'$ over $\RV$, in a way that 
yields a   reduction of the imaginaries  $Iso(\G;L)$ to imaginaries     $Iso(\G';L)$,
uniformly over   Henselian valued fields $L$ with various theories.

Note that a morphism $f: \G \to \G'$ yields, for any $L \models HF_0$, a map
$f_*: Iso(\G,L) \to Iso(\G',L)$.  We say that $f$ is an {\em elementary reduction}
of $\G$ to $\G'$ (respectively, of $\G'$ to $\G$) if $f_*$ is injective (resp., surjective.)
A {\em reduction} is a finite sequence of   elementary  reductions.

Let $G$ be a definable group acting on a definable set $V$.  Define a groupoid
$\G= \G(G,V)$ whose objects are the points of $V$.  The morphisms $v \to v'$ are defined to be:
$Mor_\G(v,v') = \{g \in G: gv = v' \}$.  Composition is multiplication in $G$.
Then $Iso(\G,L)$ is precisely the orbit space $V(L)/G(L)$.

Let $\G$ be a groupoid, with one isomorphism class.  Then all isomorphism groups
$G_a := Mor(a,a)$ are isomorphic to each other, non-canonnically: given
$a,b \Ob_\G$, choose $f \in Mor(a,b)$; then $g \mapsto f \circ g \circ f \inv$ is an isomorphism $G_a \to G_b$.   The isomorphism $G_a \to G_b$ is defined up to conjugation; if the $G_a$
are Abelian, this isomorphism is canonical, so all $G_a$ are canonically isomorphic to a 
fixed group $H$.  (This is not essential to the discussion that follows, but simplifies it.)
Let $N$ be a normal subgroup of $H$.  We   define a quotient groupoid $\G/ N$.
It has the same objects as $\G$, but the morphism set is $Mor_{\G/N} (a,b) = Mor_\G(a,b) / N$.
There is a natural morphism $\G \to \G/N$.

Our reductions will use a sequence of canonical normal subgroups of a torus $T$. 
First let $T=G_m^r$, a split torus.  The valuation map induces a homomorphism
$T \to \G^n$, with kernel $ (\Oo^*)^r$ (where $\Oo$ is the valuation ring.)
Next, we have a reduction map $(\Oo^*)^r \to (\k^*)^r$, with kernel $(1+\Mm)^r$.
Now if $T$ is any torus defined over  a valued field $F$, by definition there exists a finite Galois
extension $L$ of $F$, and an isomorphism $f: T \to G_m^r$ defined over $L$.
It is easy to see that $N:= f \inv (\Oo^*)^r$ and $N^{-1} : = f \inv (1+\Mm)^r$ do not depend
on $f$; so these are quantifier-free definable subgroups of $T$.  The quotient $T/N$ is
internal to $\G$, while $N/N^{-}$ is internal to the residue field (and $N$ is generically metastable.)   Note that $N^{-1}$ is a uniquely divisible Abelian group.  
%If $(1+m)^k=1+km +m^2(...) = 0$ then $m=0$.

We will assume $L$ is not trivially valued, so that $L^{alg} \models ACVF$.
In particular all definable torsors have $L^{alg}$-definable points.  The proof in the trivially
valued case is left to the reader.  \footnote{If   If $F$ is trivially valued, then any element of $\G$ defined over $F$ equals zero; the only   definable subgroups  of the additive group $G_a$
are thus $0,\Oo$ and $G_a$; the only definable $\Oo$-subtorsor  of $G_a$ is $\Oo$ (consider
valuation of elements).}

\<{thm}  \lbl{torus}  Let $G$ be a linear algebraic group,  $T \leq G$ a  torus.
Then $\G(G,G/T)$  reduces to a groupoid defined over $\RV$.  \>{thm}

\prf   Let $V=G/T$, $\G= \G(G,V)$.  
Consider first $\G/ N$.  Each automorphism group $Mor(a,a)$ is a uniquely divisible Abelian group
(isomorphic, over additional parameters, to $\G^r$.)  Hence given a finite subset of $Mor(a,b)$,
it is possible to take the average, obtaining a unique point. 
 In this way we can find for
each $a \in V$ a definable point $c(a) \in Mor(1,a)$ 
(where $1$ is the image in $G/T$ of $1 \in G$ .)    Given $a,b \in V$ let $c(a,b) = c(b) c(a) \inv \in Mor(a,b)$.  Then we have a subgroupoid of $\G/N$ with the same objects, and whose only
morphisms are the $c(a,b)$.

Let $\G_1$ have the same objects as $\G$, and 
$Mor_{\G_1}(a,b) = \{f \in Mor_{\G}(a,b): f N = c(a,b) \}$.  The inclusion morphism
$\G_1 \to \G$ induces  surjective maps $Iso(\G_1,L) \to Iso(\G_2,L)$ for any $L$, being surjective
on objects.  Thus $\G$ reduces to $\G_1$.

Let $\G_2 = \G_1/ {N^-}$.  The automorphism groups of $\G_2$ are isomorphic to 
$\mathfrak{t}  := N/{N^-}$, a torus over $\RV$ (i.e. a group isomorphic, with parameters,
to $(\k^*)^r$.)   

\Claim{1}  Let $L \models HF_0$.  Let $Y$ be principal homogeneous space for $N$
(defined in $ACVF_L$).  If $Y/{N^-}$ has a point in $\dcl(L)$, then so does $Y$.

\prf  Let $Y^-$ be an $L$-definable point of $Y/N^{-}$, i.e. an $L$-definable  $N^-$-subtorsor
of $Y$.  As above, since $N^-$ is uniquely divisible, $Y^-$ has an $L$-definable point. \eprf

Applying this to $Mor_{\G_1}(a,b)$,  we see that if $a,b$ are $\G_2$-isomorphic then
they are $\G_1$-isomorphic; so the natural morphism $\G_1 \to \G_2$ is injective on 
isomorphism classes over any $L \models HF_0$.   Thus $\G_1$ reduces to $\G_2$.

Finally we reduce the objects.  We have $Mor_{\G_2} (a,b) \subseteq \RV$.  
By stable embeddedness of $\RV$ there exists a definable map $j: \Ob_{\G_2} \to Y \subseteq \RV$ such that $Mor(1,a)$ is $j(a)$-definable.  It follows that
$Mor(a,b) = Mor(1,a) \times_{\mathfrak{t}} Mor(1,b)$  is $j(a),j(b)$-definable. 
Let $\G_3$ be the groupoid with objects $j(Ob_{\G_2})$, and the same morphism sets as $\G_2$.   The natural morphism $\G_2 \to \G_3$ ($j$ on objects, identity on morphisms) 
is bijective on $Iso$, since $Mor_{\G_3} (j(a),j(b)) = Mor _{\G_2}(a,b)$.  Thus $\G_2$
reduces to $\G_3$.  But $\G_3$ is over $\RV$.

\eprf

% Let $L$ be a finite Galois extension of $F$ such that $Y(L) \neq \emptyset$,
%$G= Aut(L/F)$.    We have $L \subseteq \dcl(\RV)$ (cf. e.g. \lemref{tame2} for a stronger statement.)

 Since we used only the quantifier elimination of ACVF, rather than HF, this method
 of investigation is not blocked in positive characteristic.    
  \thmref{torus} should go through for tori that split in a tamely ramified extension, 
  replacing the  unique divisibility
 argument for the additive groups by   Hilbert 90 
 for the residue field, and an appropriate extension to $\RV$.  The right statement in the general case may give a lead
with respect to  motivic integration in positive characteristic.

\>{section}

\<{section}{Galois sorts}  %\S 4

 Let $F$ be any field, and consider the 2-sorted structure $(F,F^{alg},+,\cdot)$.
 Working in the structure $(F,F^{alg})$ is convenient but harmless, since no new structure is induced on $F$ \footnote{as one easily sees by an automorphism argument.}

Let $T_0$ be a theory of   fields (possibly with additional structure.)  Let $T$ be the theory whose models have the form $(F,K,+,\cdots)$, with 
$K$ an algebraically closed field and $F$ a distinguished subfield (possibly with additional
structure) such that $F \models T$.  We can restrict attention to $K=F^{alg}$ since
$(F,F^{alg}) \prec (F,K)$.      In this section, $T$ is fixed, and ``definable'' means:  definable in $T$, imaginary sorts included.  Let $F_0 = \dcl(\emptyset)_T$.

\ssec{Definition of Galois sorts}.   Let $\mathcal{E}_n$ be the set of Galois extensions of $F$ of degree $n$, within $F^{alg}$; this is clearly 
a definable set of imaginaries of $(F,F^{alg})$.  For $e \in \mathcal{E}_n$ coding an extension $F_e$ of $F$,  let $G_e$ be the Galois group 
$Aut(F_e/F)$.  Let $\fG_n$ be the disjoint union of the $G_e$; it comes with a map $\fG_n \to \mathcal{E}_n$.
Let $\fG = (\fG_n: n \in \Nn)$.   See the Appendix for a definition at a greater level of generality,
including some definitions for Galois cohomology.

\ssec{A finiteness statement for  $H^1$ }
 
Let $A$ an algebraic group defined over $F_0$,  not necessarily commutative.   
We are interested in the first Galois
cohomology set $H^1(F,A)=H^1(Aut(F^{alg}/F,A(F^{alg})))$, where $F \models T$.  

To say that an object such as  $H^1(F,A)$ is definable means that there 
exists a definable set $H$  of $T^{eq}$ and for any $F' \models T$, a canonical
bijection $H(F') \to H^1(F',A)$.   By standard methods of saturated models, 
such a definable set $H$, if it exists, is unique up to a definable bijection. 
Given a property $P$ of definable sets (invariant under definable bijections), we 
say that $H^1(F,A)$ has $P$ if $H$ has $P$.

\<{thm}\lbl{G0}  Let $F$ be a perfect field.  If $A$ is a linear group,
 then $H^1(F,A)$ is definable, and $\fG$-analyzable.  

\>{thm}  

This will be proved as \propref{G-L7} below.

In case the Galois group of $L$  {\em property F} in the sense of \cite{serre}, (or
{\em bounded} in the sense of \cite{PAC}),   \thmref{G0} says simply
that  $H^1(F,A)$ (resp. the kernel  $H^1(F,A) \to H^1(F,G)$) is finite.  This is Theorem 5 of   \cite{serre}, Chapter 3, \S 4.  
 
 \<{question}   Is $V(K)/G(K)$ in fact internal to $\fG$?   \>{question}

Presumably it is not the case, in general, that $V(K)/G(K) \subseteq \dcl(\fG)$, even over
$\acl(0)$; 
It would be good to have an example.

 If $L$ is a Galois extension of
$F$, let $Z(L;A)$ be the set of maps $a: Aut(L/F) \to A(L)$ satisfying $a(st) = a(s) s( a(t))$.
 Two elements $a,a'$ are cohomologous if $a'(s) = b \inv a(s) s(b)$ for some 
$b \in A(L)$.
(Note $a'(1)=1=b \inv s(b)  $ so $b \in Fix(s)$ whenever $a(s)=a'(s)=1$.)  
 The quotient of $Z(L,A)$ by this equivalence relation is denoted $H^1(L/F,A)$.
  If $L \leq L'$ we have a natural map  $Z(L,A) \to Z(L',A)$, obtained by composing
 with the canonical quotient map $Aut(L'/F) \to Aut(L/F)$.  This induces a map
  $H^1(L/F,A) \to H^1(L'/F,A)$, which is injective.  

Let $Z(n;A)$ be the disjoint union of the sets $Z(L;A)$ over all Galois extensions
$L$ of $F$ with $[L:F] | n$.   Define an equivalence relation $E$ on $Z(n;A)$:
if $a_i \in Z(L_i;A)$ for $i=1,2$, write $a_1 \sim a_2$ if $a_1,a_2$ are
cohomologous in $L_1L_2$; equivalently, for some Galois extension $L$ of $F$
containing $L_1,L_2$, there exists $b \in A(L)$ such that for $s \in Aut(L/F)$, 
 $a_1(s|L_1) = b \inv a_2(s| L_2) s(b)$.
The second formulation shows that $\sim$ is an equivalence relation; the first shows
that $\sim$ is definable.  Definability of $Z(n;A)$ is clear. Denote the quotient $Z(n;a)/\sim= H(n;A)$.  

  If $n | n'$ we obtain an injective map $H(n;A) \to H(n';A)$.  It is clearly definable.  
  For any $F \models T$, $H(n;A)(F)$ is the set of elements of $H^1(F,A)$ 
  in the image of $H^1(L/F,A)$ for some Galois extension $L$ of degree $n$.
%$$H^1(F,A) = H^1(Aut(F^{alg}/F),A(F^{alg})) = \lim_{n} H(n;A)(F)$$
%and $H(n;A)(F)$ %

More generally, if $\{A_y\}$ is a definable family of algebraic groups, for $b$ from $F$
the set $H(n;A_b)(F)$ of elements of $H^1(F,A_b)$  in the image of $H^1(L/F,A)$ for some Galois extension $L$ of degree $n$ is definable uniformly in the parameter $b$.  

 \<{prop} \lbl{G-L6} Let $A$ be a linear algebraic group.  Then $H^1(F,A)$ is definable;
 for some $n$, $H^1(F,A)=H(n;A)(F)$.  
 \>{prop}
 
\prf  We have $A \leq GL_n$.  By \cite{serre-loc} Chapter X, Prop. 3, $H^1(F,GL_n)=1$.
The kernel of $H^1(F,A)  \to H^1(F,GL_n)$ is canonically isomorphic to $GL_n (F) / A(F)$.
Since $GL_n/A$ is a definable set, by a standard compactness argument, 
$\lim_n H(n;A) = \union_n H((n+1)!; A) \setminus H(n!;A)$ must also be a definable set,
i.e. for large enough $n$ the set $H((n+1)!; A) \setminus H(n!;A)$ must be empty.  \eprf

\<{cor} \lbl{G-L7f}  Let $A$ be a finite (linear algebraic) group.  Then $H^1(F,A)$
is contained in $\dcl(A,\fG)$. \>{cor}

\prf  A function from a finite set into the finite set $A$
is definable over the elements of $A$ and the elements of the domain.   
Thus $Z(n;A) \subseteq \dcl(\fG_n,A)$.  \eprf

We give a second proof, similar to the proof of (a) implies (b) in \cite{serre}, 4.1, Proposition 8.
This second proof goes through in a more general context, see the Appendix.  
 
\<{lem} \lbl{G-P7'} Let $A$ be a finite  definable group, of order $n$.  Then $H^1(F,A)=H(n!^{n!};A)(F)$.
 \>{lem}

\prf  
Let $L_0$ be the Galois extension of $F$ generated by the $n$ points of $A$.
Since by a trivial estimate $|Aut(A)| \leq (n-1)!$, we have $[L_0:F] \leq (n-1)!$.
Let $L$ be any Galois extension of $F$, containing $L_0$,   $G= Aut(L/F)$, and let
$a \in Z(L ; A)$.  We have to show that the class of $a$ in $H^1(L/F,A)$
is in the image of some $H^1(L'/F,A)$ with $[L':F] \leq (n!)^{n!}$.  In fact we will prove
this even at the level of cocycles.  The restriction of $a$ to 
$G_0=Aut(L/L_0)$ is a homomorphism $a | G_0: G_0 \to A$.   Let $G_1$ be the kernel
of $a | G_0$.  Then $[G_0:G_1] \leq n$, so $[G:G_1] \leq n!$.  The number of conjugates
of $G_1$ in $G$ is thus $\leq n!$; their intersection $G_2$ has index $\leq (n!)^{n!}$.  
Let $L'$ be the fixed field of $G_2$.  Then $a$ factors through a function $a'$ on
$G/G_2$, so $a' \in Z(L'/F,A)$, as required.  \eprf 
 
\<{cor} \lbl{G-L7T}  Let $A$ be a torus, i.e. $A$ becomes isomorphic to $G_m^n$
after base change to some Galois extension $F'$.  Say $[F':F]=m$, and
let $B = \{a \in A: a^m =1 \}$.   Then $H^1(F,A) \subseteq  \dcl(B,\fG)$.  If $F$ is perfect, then 
this holds for any connected solvable $A$ (for an appropriate finite group $B$.)
\>{cor}

\prf The proof in \cite{serre}, Chapter III, Theorem 4, goes through, showing
that $H^1(F,A) \cong H^1(F,B)$.  \eprf

\<{rem} \lbl{G-L7TR} We also have $H^1(F,A) \subseteq \dcl(F)$, since the Galois group $Aut(F^{alg}/F)$ 
acts trivially on $H^1(F,A)$.  \>{rem}

\<{prop}  \lbl{G-L7}  Assume $F$ is perfect.  Let $A$ be a linear algebraic group.  Then $H^1(F,A)$ is $\fG$-analyzable.  \>{prop}

\prf  The proof of \cite{serre}, chapter III, Theorem 4 goes through.   \eprf

\<{rem}  The proof of \propref{G-L6} shows that for any embedding of algebraic groups $A \to B$, the kernel of $H^1(F,A) \to H^1(F,B)$ is definable.  It can be shown as in \cite{serre}, chapter III, Theorem 5 that over a perfect field, this kernel is     $\fG$-analyzable.  \>{rem}

%E.g. by Lemma 3.5 of \cite{H-B},  the family of all subgroups $S$ of a linear group $G$
%equal to the connected component of their normalizer, is a uniformly definable family.
%In particular the Cartan subgroups of $G$ - but all eventually conjugate, anyhow.

\<{cor} \lbl{G2}   Let $L$ be a perfect field.  
Let $G$ be a linear algebraic group over $L$, and let $V$ be a homogeneous
space for $V$ defined over $L$, i.e. $G$ acts on the variety $V$ and $G(L^{alg})$
acts transitively on $V(L^{alg})$.  
Then in the structure $(L,L^{alg},+,\cdot)$, 
 $V(L)/G(L)$ is ${\fG}_L)$-analyzable.  \>{cor}

\prf  Immediate from \propref{G-L7}, since after picking any point  $c \in V(L)$ and letting
$H = \{g \in G: gc=c \}$, 
we have a $c$-definable injective map $V(L)/G(L) \to H^1(F,H)$.  \eprf

\ssec{Henselian fields}

We now move to valued fields of residue characteristic zero.

\<{lem} \lbl{tame2} Let $K$ be  a Henselian field with residue field of characteristic $0$.   Let $\k$ denote the residue
field, $\mu=\union_n \mu_n$ the roots of unity in $\k^{alg}$, $\G$ he value group.  
Then in $(K,K^{alg})$
we have $\fG_K \subseteq \dcl( \fG_\k \union \mu \union \union_{n} \G/ n \G  )$
\>{lem}  

\prf    Let $L$ be a finite Galois extension of $K$.  We have to show that $Aut(L/K)  \subseteq \dcl( \fG_\k, \mu, \G/ n \G )$ for some $n$.

We will use some standard valuation theory.  Call $K'$ a {\em ramified root extension} of $K$ if it is a finite
purely ramified extension obtained by adding roots to some elements of $K$.

\Claim{1} There exists a finite unramifield extension $L'$
of $L$, 
 an unramified Galois extension $K_{u} \leq L'$ of $K$,
 and a ramified root extension $K_{rr} \leq L'$ of $K$, such that
%$L'$ is a ramified root extension of $K_u$,  
$ L'/K_{rr}$ is unramified, and
$Aut(L'/K) = Aut(L'/K_u) Aut(L'/K_{rr})$.

\prf The finite group $\G(L)/\G(K)$ is a direct sum of finite cyclic groups, $\oplus_{i=1}^k \Zz / n_i \Zz$; let $c_1,\ldots,c_k \in L$
be such that $\val(c_1)+ \G(K), \ldots, \val(c_k) + \G(K)$ are generators for these cyclic groups.
So $\val(c_i^{n_i} ) = \val(d_i)$ for some $d_i \in K$.  Let $e_i =  c_i^{n_i } / d_i$, so that $\val(e_i)=0$.  In a Henselian field of residue characteristic prime to $n$, if an element $f$ has $\val(f)=0$, and $\res(f)$ has an $n$'th root, then by Hensel's lemma so does
$f$.  Hence 
in some unramifield Galois extension $L'$ of $L$, each $e_i$ has an $n_i$'th root.  
Clearly $d_i$ has an $n_i$'th root $r_i \in L'$.  Let $K_{rr}= K(r_1,\ldots,r_k)$.  
 Let $K_u$ be the maximal
unramified subextension of $L'$.  Then $L'$ has the same residue field over $K_u$;
$L'/K_u$ is purely ramified.  Since
$\val(L')=\val(L)= \val(K_{rr})$, $L'/K_{rr}$ is unramified.  In particular
$K_u \meet K_{rr} = K$, so $Aut(L'/K) = Aut(L'/K_u) Aut(L'/K_{rr})$ by Galois theory. \eprf

\Claim{2}  $ Aut(L'/K_{rr})  \subseteq \dcl( \fG_\k)$.

\prf  The canonical homomorphism $Aut(L'/K_{rr}) \to Aut(\k_{L'} / \k)$ is an isomorphism,
since $L'/K_{rr}$ is unramified.  This homomorphism is definable, hence embeds
the elements of $Aut(L'/K_{rr})$ into $\fG_\k$. \eprf

\Claim{3}    $ Aut(L'/K_u) \cong Hom(\G(L)/\G(K), \mu_n)$ (canonically and definably.)

\prf  Let $E = \{c \in L^*: c^n \in K^* \}$.  
 Define a map $b: Aut(L' / K_u) \times E \to \mu_n$ by $b(\si,e) = \si(e)/e$.  
Clearly $b$ is multplicative in the second variable.
%Let $\Oo^*=\{x: \val(x)=0 \}$.  
% We have a residue homomorphism $E \meet \Oo^* \to \res(L')=\res(K_u
 If $e \in E$
and $\val(e)=0$, then $\res(e^n)$ has an $n$'th root in $\res(L')$ and hence
in $\res (K_u)$, since $L'/K_u$ is purely ramified.  So $e^n$ has an $n$'th root in $K_u$.
Since all roots of unity in $L'$ lie in $K_u$, we have $e \in K_u$; hence
$b(\si,e)=1$ for all $\si$.   
More generally if  $e \in E$ and $\val(e) \in \G(K)$, 
then $\val(e/c)=0$ for some $c \in K^*$, so $b(\si,e/c)=0$ and hence $b(\si,e)=0$
for all $\si$.  
Thus $b$ factors through $b': Aut(L'/K_u) \times  \G(L')/\G(K) \to \mu_n$.
We obtain a homomorphism $Aut(L'/K_u) \to Hom(\G(L')/\G(K) , \mu_n)$.
It is injective since if $b(\si,r_i)=1$ for each $i$, then $\si(r_i)=r_i$ so $\si=1$.  
Surjectivity comes from cardinality considerations (but will not be needed.) 
\eprf

\Claim{4}  $ Aut(L'/K_u) \subseteq \dcl( \mu_n  , \G/ n \G )$ where $n=[L':K_u]$.

\prf  By Claim 3, $Aut(L'/K_u)$ is definably isomorphic to  
$Hom(\G(L)/\G(K), \mu_n)$.  Now $\G(L)/ \G(K)$ is a finite subgroup of $(1/n) \G(K) / \G(K)$,
hence isomorphic to a finite subgroup   $S \leq \G(K) / n \G(K)$.  Each element
of $S$ lies in $\G / n \G$, and so a homomorphism $S \to \mu_n$ is definable
over teh elements of $S$ and the elements of $\mu_n$, and the claim is proved.  \eprf

The lemma follows from Claims 1,2 and 4.
\eprf

\<{thm}   \lbl{CDm} Let $K$ be a Henselian field with residue field of characteristic $0$.  Let $G$ be an algebraic group over $K$, acting on a variety $V$; assume $G(K^{alg})$ acts transitively on $V(K^{alg})$.  Let $V(K)/G(K)$
be the orbit space.   Then for some $n$, $V(K)/G(K)$ is analyzable over  sorts $\G/n\G$ and
the $n$-th Galois sort   $\fG_n(\k)$ of the residue field.
\>{thm}

\prf[Proof of \thmref{CDm}]  Given $c \in V(K)$, $H=\{g: gc=c\}$, there is a canonical $c$-definable bijection between $V(K) / G(K)$
and the kernel of $H^1(K,H) \to H^1(K,G)$.  
Hence the theorem follows from \thmref{G0} and   \lemref{tame2}.
\eprf

\<{rem}  
If we add to the structure the Denef-Pas splitting, 
(at least when e.g. $\G$ is a $\Zz$-group, but probably in general),  
it follows that $V(K)/G(K)$ is tame,
i.e. $V(K)/G(K) \subseteq \dcl(\k,\G/n\G)$.     \rm Indeed for every imaginary definable set $D_0$ of $\G$, 
either $D_0 \subseteq \dcl(\G/ n \G)$ for some $n$, or else $D \subseteq \dcl(D_0)$
for some infinite definable $D \subseteq \G$.  But it is easy to see that $D$ is not
internal over $\k \union \G/m\G$.  Hence the definable sets internal over
 $\k \union \G/m\G$ are contained in $ \dcl(\k \union \G/ n \G)$ for some $n$, and by
 induction the same goes for analyzability.   
 \>{rem}

%\<{rem}  \rm We presented three different descriptions of the orbit space $V/G$
%(finite imaginaries, groupoids, Galois sorts.)
%While our emphasis is on the methods, we compare here the strength of the statements 
%obtained with respect to coding the orbits in $\RV$.  
%The proof of \S1-3, using finite imaginaries, applies only when $V(L)/G(L)$ is finite.
%This restriction is not needed for the  proof in \S 4 and \S 5.  In \S 4 we  $V=G/T$ for some torus $T$; and we show only
%that the orbit space is coded in $\RV$.   In \S 5 we do not assume
%this.
%None of the proofs assume $V$ is defined over a number field, or that the residue field is 
%pseudo-finite.  
%  \>{rem}  

 \>{section}

\<{section}{Appendix}  \lbl{appendix}

It may seem at first sight that Galois sorts are peculiar to fields; but in fact
they can be defined for any theory eliminating imaginaries on its finite subsets, see below.  

In particular, the Galois sorts of $\RV$ will be defined.  This will permit 
the    remark that, for a theory $T$ of Henselian fields of residue characteristic $0$,
if $T_{\RV}$ is the induced theory on $\RV$, we have 
${\fG}_T \cong {\fG}_{T_{\RV}}$; which, along with \thmref{G0}, clarifies the reductions of \cite{CD} and the present paper.  

Fix a language $L$.  
Let $\tT$ be a theory admitting elimination of imaginaries with respect to finite sets of tuples. 
Thus for    any   finite product $S$ of sorts of $L$, and any $m \in \Nn$, 
we have given a sort  $S[\leq m]$ and a  function $c_{S,m}: S^m \to S[\leq m]$, 
such that $$\tT \models (c_{S,m}(x) = c_{S,m}(y) \iff \bigwedge_{\si \in Sym(m)} x^\si = y)$$
Let $S[m]$ be the image of the distinct $m$-tuples.  

For simplicity we assume also that $\tT$ eliminates quantifiers, and that 
any quantifier-free definable function of $\tT$ is given by a term (piecewise), so 
that substructures are definably closed.    Let $T_0$ be the theory of  substructures of models of $\tT$.  If 
$M \models T_0$, let $M^{alg}$  denote the algebraic closure of $M$ within some
 model of $\tT$.  The Galois imaginaries, strictly speaking, belong to $T_0^a = Th(\{(M,M^{alg}): M \models T_0 \})$.  We will also note some related imaginary sorts of $T_0$ itself.
 In practice it seems more convenient to use $T_0^a$, then note by considerations of 
 stable embeddedness that  definable sets on which the Galois group acts trivially
 belong to $T_0$.  For instance this is the case with the cohomology sets $H^1(Aut(M^{alg}/M), A)$ considered below.

Below we abuse notation, identifying elements of $S[m]$ with $m$-element sets. 
Note however that if $M \models T$, then 
$S[m](M)$ corresponds to $m$-element subsets of $M^{alg}$, not of $M$.

We define the {\em Galois sorts } $\fG$ 
of $T_0$.  These are certain definable sets of imaginaries (not quantifier-free.)

   We have a relation $(\subseteq) \subset S[\leq m]^2$,
corresponding to inclusion of finite sets.
 Let $S_{irr}[m]$ be the set of {\em irreducible} elements of $S[m]$:
 $$x \in S_{irr}[m] \iff ( x \in S[m] \wedge \bigwedge_{1 \leq k < m} \neg ( \exists y \in S[k]) (y \subseteq x))$$
Thus if $M \models T_0$,  then $S_{irr}[m](M)$ corresponds to the set of orbits of $Aut(M^{alg}/M)$
on $M^{alg}$ of size $m$.  
If $k \leq m$, for $x \in S_{irr}[k], y \in S_{irr}[m]$ we let $Mor(x,y)$ be the set of codes of  functions $y \to x$ (viewed as subsets of $y \times x$.)   This makes $\union_{S,m} S_{irr}[m]$ into a category $\fC$.  

If $a \in S(M^{alg})$ has $m$ conjugates under $Aut(M^{alg}/M)$, let $s(a)$ be the
(code for) the set of conjugates.  Then $s(a) \in S_{irr}[m](M)$, and every element
of $S_{irr}[m](M)$ arises in this way.  
 A $\fC$-morphism $s(b) \to s(a)$  exists iff $M(a) \leq M(b)$.  In particular, 
 $s(a),s(b)$ are $\fC$-isomorphic iff    $M(a)=M(b)$, i.e. $a,b$ generate the same substructure of $M^{alg}$ over $M$.  
Isomorphism classes of $\fC[M]$ correspond to 
finitely generated
extensions of $M$ within $M^{alg}$.  

 If $M \models T_0$ and $a \in S_{irr}[k](M), b \in S_{irr}[m](M)$, then $Mor(a,b)(M)$ is a finite set, possibly empty.  By irreducibility, any function coded 
in $Mor(a,b)[M]$ must be surjective.   In particular if $k=m$ it must be bijective, so the 
subcategory with objects $S_{irr}[m](M)$ is a groupoid (every morphism is invertible.)
For $s \in S_{irr}[m]$, let $H_s = Mor(s,s)$.  
 Let $S_{gal}[m]$ be set of $s \in S_{irr}[m]$ such that $H_s$ acts regularly on $s$.
   Let $\fE_{S,m}$ be the set of $\fC$- isomorphism classes of 
objects in $S_{gal}[m]$.  If $s \in S_{gal}[m]$ codes a set $D_s$, let $Gal(s)=Aut(D_s;H_s)$
be the opposite group to $H_s$, i.e. the group  of permutations of $D_s$ commuting
with each element of $H_s$.  

Let $E_s = Gal(s)^m / ad$ be the set of $Gal(s)$-conjugacy
classes on $Gal(s)^m$.  (More canonically, we could look at $E^k_s = Gal(s)^k / ad$
for all $k \in \Nn$, 
but   all $E^k_s$ are definable from $E^m_s = E_s$.)   

The definition of $Gal(s)$ requires knowledge of $D_s$, 
i.e. of a particular choice of an algebraic closure $M^{alg}$ of $M$.  But $E_s$ is canonical
and does not depend on this choice.  

 If $s,s'$ are isomorphic, i.e. $Mor(s,s') \neq \emptyset$,
the choice of $f \in Mor(s,s')$ yields an isomorphism $Gal(s) \to Gal(s')$
(namely $\tau \mapsto f \circ \tau \circ f \inv$), and this isomorphism does not depend
on the choice of $f$.  Indeed let $M(s)$ be the substructure of $M^{alg}$ generated over
$M$ by any element of $s$;  the substructures $M(s)$,$M(s')$ are equal, and
$Gal(s) = Aut(M(s)/M)$.   
  The induced bijection $E_s \to E_{s'}$ depends therefore neither
on this choice nor on a choice of $M^{alg}$.  

Let ${\equ{\fC}}$ denote $\fC$- isomorphism, and let  $\fE_{S,m} =  S_{gal}[m] / {\equ{\fC}}$.
 then if    $e = s /  {\equ{\fC}} \in \fG_{S,m}$ we
may define $Gal(e) = Gal(s) , M(e)=M(s), E_e = E_s$.  Let $\fG_{S,m} $ be the disjoint union over
$e \in \fE_{S,m}$ of the groups $Gal(e)$, and let $\fG_{S,m} \to \fE_{S,m}$
be the natural map.   Similarly let $\cG$ be 
 the direct limit over all $S,m$ of $\fE_{S,m}$.  

Let $\cG$ denote the family of all sorts $\cG_{S,m}$ and $\fE_{S,m}$.   These will be  called
the {\em Galois sorts of $T_0$}. 
Let $T_0^{gal}$ consist of all sorts interpretable over these sorts.  (I.e. close under
quotients by definable equivalence relations.)   Note that   the groups $Gal(e)$ 
are not themselves part of $\cG$, in general.

 \<{rem} 
The Galois
group $G=Aut(M^{alg}/M)$ is the projective limit over $\fE_{S,m}$ of all groups $Gal(e)$.  
An element
$e \in \fG_{S,m}(M)$ corresponds to a   normal open subgroups $N(e)$ of $G$, and
we have $Gal(e) \cong G/N(e)$ canonically.  
 Similarly for any $k$ the space $G^k / ad_G$ of conjugacy classes of $G$
on $G^k$ can be deduced from $\cG_{S,m}$ by an appropriate projective limit.
  \>{rem}

\<{rem}
Let $M \models T$.  Let $G = Aut(M^{alg}/ M)$, $C$ the centralizer of $G$ in $Aut(M^{alg})$.
Then  $Aut(M / \cG(M))$ is the image of $C$ in $Aut(M)$.     Equivalently, 
$\si \in Aut(M / \cG(M))$ iff every extension  $\tau$ of $\si$ to $Aut(M^{alg})$ 
normalizes $Aut(M^{alg}/M)$, and induces an inner automorphism of this group.
 It follows  in particular that any such $\tau$
preserves each Galois extension of $M$.  \>{rem}

\<{rem} \lbl{PE} Assume (PE):   $L$ has a  sort $S_{basic}$ for which the primitive element theorem is
valid, i.e. any finite extension of $M \models T_0$ is  generated by a single element of $S_{basic}$.  Then
only the Galois sorts $\fE_m= \fE_{S_{basic},m}$ and $\cG_{m}= \fE_{S_{basic},m}$
need be considered; any Galois sort is in definable bijection with a definable subset of these. \>{rem}

\ssec{Galois cohomology}

Let  $G$ be a profinite group, acting continuously on a discrete group $A$.  We write $^g a$
for the action.  

Recall  the definition of $H^1(G,A)$ (\cite{serre} Chap. I, 5.1.)  A 1-cocyle
is a continuous map $a: G \to A$ satisfying $a(st) = a(s) ^s a(t)$.  The set of 1-cocycles is denoted
$Z^1(G,A)$.  Two cocycles $a,a'$ are cohomologous if $a'(s) = b \inv a(s) ^s b$ for some 
$b \in A$.   The quotient of $Z^1(G,A)$ by this equivalence relation is denoted $H^1(G,A)$.
 The action of $G$ on $Z^1(G,A)$, induced
from the actions on $G$ and on $A$, satisfies:  $^g a (t) = b \inv a(t)  ^t b$, where $b = a(g)$;
hence $^g a$ is cohomologous to $a$, so $G$ acts trivially on $H^1(G,A)$.  

When $A$ is a definable group of   $\tT$, for any 
$M \models T_0$ we have a profinite group $G_M:= Aut(M^{alg}/M)$, and a
continuous action of $G_M$ on $A(M^{alg})$.  We write $G$ for the functor $M \mapsto G_M$.

\<{lem}  Assume (PE) (cf. \ref{PE}).  Let $A$ be a definable finite group of $\tT$.  Let $G_M$ denote the automorphism
group $Aut(M^{alg}/M$.    Then $H^1(G, A)$ is  interpretable in the Galois sorts of $T_0$.
In other words there exists a definable set $S$ of $T_0^{gal}$ and for any $M \models T_0$
a canonical bijection $S(M) \to  H^1(G_M,A)$.  \>{lem}

\prf  The proof of \lemref{G-P7'} goes through. 
\eprf

\>{section}

\>{document}

\>{document}